\documentclass[final]{siamltex}     
\usepackage{epsfig}
\usepackage{amsfonts}
\usepackage{ifpdf}
\usepackage[backref,citecolor=blue,%
    pdftitle={Bounds on the convergence of Ritz values from Krylov %
      Subspaces to interior eigenvalues of Hermitean matrices},%
    pdfauthor={Chris Johnson and A. D. Kennedy},%
    pdfkeywords={Lanczos, Krylov, spectrum, eigensolvers, linear algebra},%
    pdfstartview={FitH},%
    pdfpagemode={FullScreen}]{hyperref}
\ifpdf
\let\mysection=\section
\def\section#1{\mysection{#1}\pdfbookmark[-1]{#1}{\thesection}}
\let\mysubsection=\subsection
\def\subsection#1{\mysubsection{#1}%
  \pdfbookmark[0]{#1}{\thesection.\thesubsection}}
\let\myappendix=\appendix
\def\appendix{\myappendix\pdfbookmark[-1]{Appendix}{\thesection}}
\fi
\def\epspdffile#1{\leavevmode\ifpdf\epsffile{#1.pdf}\else\epsffile{#1.eps}\fi}
\let\definition=\emph               
\def\Lanczos{L\'anczos}             
\def\Hermitean{Hermitean}           
\def\Span(#1){\mathop{\rm span}\nolimits(#1)}   
\def\Krylov(#1,#2,#3){{\mathcal K}_{#3}(#1,#2)}	
\def\binomial(#1,#2){\left(\begin{array}{c} #1 \\ #2 \end{array}\right)}
\def\implies{\Rightarrow}	    
\def\intersect{\cap}		    
\def\defn{\equiv}		    
\def\mtx#1{{\mathsf #1}}            
\def\vect#1{{\mathbf #1}}           
\def\u{\vect u}			    
\def\v{\vect v}			    
\def\s{\vect s}			    
\def\w{\vect w}
\def\y{\vect y}			    
\def\z{\vect z}			    
\def\A{\mtx A}			    
\def\H{\mtx H}			    
\def\Q{\mtx Q}			    
\def\T{\mtx T}			    
\def\U{\mtx U}			    
\def\R{{\mathbb R}}		    
\def\rational#1#2{{\mathchoice{\textstyle{#1\over#2}}%
  {\scriptstyle{#1\over#2}}{\scriptscriptstyle{#1\over#2}}{#1/#2}}}
\def\half{\rational12}		    
\begin{document}
%
%
\title{Bounds on the convergence of Ritz values from Krylov Subspaces to
  interior eigenvalues of Hermitean matrices}
\author{Chris Johnson
  \hskip1em and\hskip1em
  A.~D.~Kennedy \\
  Tait Institute, SUPA, NAIS and EPCC, \\
  Department of Physics and Astronomy, \\
  The University of Edinburgh, The King's Buildings, \\
  Edinburgh, EH9~3JZ, Scotland}
\date{\small{\it Version of \today}}
\maketitle
\begin{abstract}
  \noindent We consider bounds on the convergence of Ritz values from a
  sequence of Krylov subspaces to interior eigenvalues of
  \Hermitean\ matrices.  These bounds are useful in regions of low spectral
  density, for example near voids in the spectrum, as is required in many
  applications.  Our bounds are obtained by considering the usual
  Kaniel--Paige--Saad formalism applied to the shifted and squared matrix.
\end{abstract}
\pagestyle{myheadings}
\thispagestyle{plain}
\markboth{C. Johnson and A. D. Kennedy}{Convergence of interior Ritz values in
Krylov subspaces}
\section{Introduction}

In this paper we derive a bound on how well an eigenvalue of a
\Hermitean\ matrix is approximated by the corresponding Ritz value obtained
from its restriction to a Krylov subspace.  Such bounds are useful in practice
because no Krylov space based method for determining eigenpairs (such as
Arnoldi or \Lanczos\ algorithms) can do better than this even in exact
arithmetic.

We shall show that the eigenvalues \(\lambda_\alpha\) of an \(N\times N\)
\Hermitean~matrix \(\A\) are approximated by its Ritz values \(\theta_\alpha\)
from the Krylov space \(\Krylov(\A,\v,n)\) with an error bounded~by
\[
  0 \leq \lambda_\alpha-\theta_\alpha
  \leq \min_\Lambda\frac{2(\lambda'_j-\lambda'_N)(K'_j\tan\angle'_j)^2
    e^{-2(n-2j)\mu'_j+O({\mu'_j}^2)}}{|\lambda_\alpha-\Lambda|}
\]
where \(\lambda'_1\geq\lambda'_2\geq\cdots\geq\lambda'_N\) are the eigenvalues
of \(\A'\defn-(\Lambda-\A)^2\) with \(j\) such that \(\lambda'_j=
-(\Lambda-\lambda_\alpha)^2\),
\[
  K'_j \approx \prod_{i=1}^{j-1}\frac{\lambda'_i-\lambda'_N}
	    {\lambda'_i-\lambda'_j},\qquad
  \mu'_j=\sqrt{\frac{\lambda'_j-\lambda'_{j+1}} {\lambda'_{j+1}-\lambda'_N}},
\]
and \(\angle'_j\) is the angle between the projection of \(\v\) onto
\(\Span(\z_{j+1},\ldots,\z_N)\) and \(\z_j\) where \(\z_i\) is an eigenvector
of \(\A'\) (and thus also of \(\A\)) belonging to the
eigenvalue~\(\lambda'_j\).

In section \S\ref{sec:applications} we briefly indicate why the problem of
finding approximations to interior eigenvalues of large \Hermitean\ matrices
is important in a variety of physical applications.  In
\S\ref{sec:basic-properties} we review some basic properties of
\Hermitean\ matrices introducing the notation to be used throughout the paper,
and in \S\ref{sec:subspace-bounds} we review some useful bounds on Ritz values
for general subspaces.  Section \S\ref{sec:KPS} derives the bounds on the
convergence of extremal eigenvalues due to Kaniel~\cite{Kaniel:1966},
Paige~\cite{Paige:1971} and Saad~\cite{Saad:1980}.  We then present a
straightforward but apparently new generalization that gives bounds for the
interior eigenvalues close to ``voids'' in the spectrum in
\S\ref{sec:s-s-bounds}, followed by some examples illustrating the application
of our bounds.

\section{Applications} \label{sec:applications}

The problem of finding good approximations to interior eigenvalues in
low-density regions of the spectrum of a \Hermitean\ matrix is common to many
areas of computational science, but the particular application that led us to
consider this problem \cite{Kennedy:2006ax} is that of evaluating the
Neuberger operator for lattice QCD (Quantum Chromodynamics being the quantum
field theory of the strong nuclear force).  This requires us to evaluate the
\definition{signum} function of the ``\Hermitean\ Dirac operator'' \(\gamma_5
D\), which is defined by diagonalizing this matrix and taking the sign
(\(\pm1\)) of each of its eigenvalues.  It is far too expensive to carry out
the full diagonalization, so we use a Zolotarev rational approximation for the
signum function as this can be evaluated just using matrix addition,
multiplication, and inversion (using a multi-shift solver for its stable
partial fraction expansion).  This approximation is expensive for very small
eigenvalues of \(\gamma_5 D\), and as there are only a relatively small number
of these we want to deflate them and take their sign explicitly.  While the
Neuberger operator is \Hermitean, it is not positive; furthermore, for
physically interesting parameters its spectrum has a large void around zero,
with perhaps only a few ``topological'' or ``defect'' modes leading to
discrete eigenvalues within this void.  The relevant numerical problem is
therefore to find the eigenpairs within or close to this void. 

Previously \cite{Edwards:2004sx} such eigenpairs had been found by finding the
extremal eigenpairs of the square of the Neuberger operator, but this then
required further effort to separate artificial degeneracies of the squared
matrix arising from (almost) degenerate positive and negative eigenvalues of
the original matrix.  We observed that there is in fact no need to square the
matrix, as the desired eigenpairs converge rapidly for the Neuberger operator
itself, despite the fact that the eigenvalues are far from either end of the
spectrum and thus the usual Kaniel--Paige--Saad bounds are not applicable.  We
used a new algorithm based on the \Lanczos\ algorithm with selective
orthogonalization (LANSO) method of Parlett and Scott \cite{Parlett:1979:LAS}
(which avoids the appearance of duplicates of eigenvalues due to rounding
errors in finite-precision arithmetic computations) that allows us to find
eigenpairs with eigenvalues of a \Hermitean\ matrix \(\A\) that lie within a
specified interval.  We shall describe this algorithm in a forthcoming
publication.

The fact that such interior eigenpairs converged rapidly led us to consider
how to find realistic bounds in this case, and we were led to use the
Kaniel--Paige--Saad bounds for the extremal eigenvalues of the shifted and
squared matrix \(-(\Lambda-\A)^2\) to give bounds on the convergence of the
original matrix \(\A\).

We do not believe that applications in which such interior eigenvalues in
low-density spectral regions are important are unusual.  To give just one more
example, we may consider the Hamiltonian for an electron in a metal,
semiconductor, or insulator, which may be approximated by a large
\Hermitean\ matrix with ``conduction bands'' separated by large ``band gaps''.
The details of the eigenvalues deep within the bands are not too important:
for example their separation is related to the inverse (macroscopic) volume of
the system under consideration, but the eigenpairs at the edges of and between
the bands determine many physically important properties.

\section{Basic Properties of \Hermitean\ Matrices}
\label{sec:basic-properties}

In this section we summarize the basic properties of \Hermitean\ (symmetric)
matrices, with the goal of specifying our notation.

A matrix \(\A\) is \definition{Hermitean} (with respect to a sesquilinear
inner product) if \(\A=\A^\dagger\), which means \((\u,\A\v)=(\A^\dagger\u,\v)
=(\A\u,\v)=(\v,\A\u)^{*},\) or equivalently \(\u^\dagger\cdot\A\v=(\A^\dagger
\u)^\dagger\cdot\v=(\A\u)^\dagger\cdot\v=(\v^\dagger\cdot\A\u)^{*}\).  An
eigenvalue \(\lambda\) of \(\A\) satisfies \(\A\z=\lambda\z\) where
\(\z\neq0\) is the corresponding eigenvector.  The eigenvalues are real
\[
  \lambda = \frac{(\z,\A \z)}{(\z,\z)}
    = \frac{(\z,\A^\dagger \z)}{(\z,\z)}
    = \frac{(\z,\A\z)^{*}}{(\z,\z)^{*}}
    = \lambda^{*},
\]
and the eigenvectors corresponding to different eigenvalues are orthogonal:
\(\lambda(\z',\z)=(\z',\A\z)=(\A\z,\z')^{*}=(\z,\A\z')^{*}={\lambda'}^{*}
(\z,\z')^{*}=\lambda'(\z',\z),\) hence \((\lambda'-\lambda) (\z',\z)=0\), so
\(\lambda\neq\lambda'\implies(\z',\z)=0\).  We can choose eigenvectors
belonging to the same eigenvalue to be orthonormal, for example by using the
Gram--Schmidt procedure, as any linear combination of such eigenvectors is
also an eigenvector.

Since any matrix can be reduced to triangular form \(\T\) by a unitary
(orthogonal) transformation\footnote{This is ``Schur normal form'', which
follows from the Cayley--Hamilton theorem that every matrix satisfies its
characteristic equation, and the fundamental theorem of algebra which states
that the characteristic polynomial \(p(\lambda) = \det(\A-\lambda)\) has
exactly \(\dim\A\) complex roots, counting multiplicity.} (change of basis),
\(\A=\U\T\U^{-1}=\U\T\U^\dagger\), and \(\T^\dagger=(\U^\dagger\A\U)^\dagger=
\U^\dagger\A^\dagger\U=\U^\dagger\A\U=\T\), it follows that \(\T\) is real and
diagonal, and thus \(\A\U=\U\T\) so the columns of \(\U\) furnish the
orthonormal eigenvectors.

\section{Subspace Bounds} \label{sec:subspace-bounds}

\subsection{Rayleigh Quotient}

We consider how the eigenvalues and eigenvectors of an approximation to a
matrix \(\A\) are related to the actual eigenpairs.  Consider the
\definition{Rayleigh quotient} which has the spectral representation
\[
  \rho(\u,\A) = \frac{(\u,\A \u)}{(\u,\u)}
    = \sum_j\lambda_j\frac{|(\z_j, \u)|^2}{\sum_k|(\z_k,\u)|^2}.
\]
This is a weighted mean of eigenvalues, so it satisfies the inequalities
\(\lambda_1\geq\rho(\u,\A)\geq\lambda_N\), where \(N=\dim\A\) and we have
ordered the eigenvalues such that \(\lambda_1\geq\lambda_2\geq\cdots\geq
\lambda_N\).

\subsection{Ritz Pairs}

The extrema of the Rayleigh quotient may be found by varying the components
\(u_i=(\z_i,\u)\) of \(\u\) subject to the constraint that \(\|\u\|=1\).
Introducing a Lagrange multiplier \(\mu\) we find \((\A-\mu)\u=0\), whose
solutions are the eigenpairs \(\mu=\lambda_j\) with \(\u=\z_j\) for
\(j=1,\ldots,N\).  If we restrict \(\u\) to an \(m\)-dimensional subspace
\(X\subseteq\R^N\) then the extrema of the Rayleigh quotient are the
\definition{Ritz pairs} \(\mu=\theta_j\) with \(\u=\y_j\) for
\(j=1,\ldots,m\).  The Ritz pairs are the eigenvalues and eigenvectors of the
restriction of \(\A\) to the subspace, \(\H=\Q^\dagger\A\Q\), that is
\(\H\s_j=\theta_j\s_j\).  Here \(\Q\) is an \(N\times m\) matrix with
orthonormal columns (\(\Q^\dagger\Q=1\)) and \(\Q\Q^\dagger\) is an orthogonal
projector onto~\(X\).  Note that the vectors \(\y_j=\Q\s_j\in\R^N\) are not
eigenvectors of \(\A\), although \(\Q^\dagger(\A-\theta_j)\y_j=0\), since the
\definition{residual} \(\mtx R=\A\Q-\Q\H\neq0\) in general.  The Rayleigh
quotient \(\rho(\u\in X,\A)\) restricted to the subspace \(X\) is a weighted
mean of Ritz values so \(\theta_1\ge\rho(\u\in X,\A)\geq\theta_m\) where we
have ordered the Ritz values \(\theta_1\geq\theta_2\geq\cdots\geq\theta_m\).
For large matrices our strategy is to compute the Ritz pairs of \(\A\) for
suitably chosen subspaces \(X\); so we are interested in how well these may
approximate the eigenvalues of \(\A\) itself.

\subsection{MaxMin and MinMax Bounds}

Let \(S_j\subseteq X\) be an arbitrary subspace of dimension \(j\), and
\(C_{j-1}\subseteq X\) another arbitrary subspace of codimension \(j-1\),
i.e., a subspace of dimension \(m-j+1\) if \(m<\infty\).  There must be a
non-zero vector \(\v\in S_j\intersect C_{j-1}\), where of course \(\v\)
depends upon \(S_j\) and \(C_{j-1}\); therefore \(\min_{\u\in S_j}\rho(\u,\A)
\leq\rho(\v,\A)\leq \max_{\u\in C_{j-1}} \rho(\u, \A)\).  Hence, for any pair
of subspaces \(S_j\) and \(C_{j-1}\), \( \min_{\u\in S_j}\rho(\u,\A)\leq
\max_{\u\in C_{j-1}}\rho(\u,\A)\).  Taking the maximum over all \(S_j\) and
the minimum over all \(C_{j-1}\) we obtain \(\max_{S_j\subseteq X} \min_{\u\in
S_j}\rho(\u,\A)\leq\min_{C_{j-1}\subseteq X}\max_{\u\in C_{j-1}}\rho(\u,\A)\).
For the particular subspace \(S_j=\Span(\y_1,\ldots,\y_j)\) the minimum of the
Rayleigh quotient is \(\theta_j\), likewise for \(C_{j-1}=\Span(\y_j,\ldots,
y_m)\) its maximum is also \(\theta_j\), so \(\theta_j\leq \max_{S_j\subseteq
X}\min_{\u\in S_j}\rho(\u,\A)\leq\min_{C_{j-1}\subseteq X}\max_{\u\in C_{j-1}}
\rho(\u, A)\leq\theta_j\):we have thus established the \definition{MaxMin} and
\definition{MinMax} bounds
\[
  \max_{S_j\subseteq X}\min_{\u\in S_j}\rho(\u,\A) = \theta_j 
    =\!\!\min_{C_{j-1}\subseteq X}\max_{\u\in C_{j-1}}\rho(\u,\A).
\]

\subsection{Krylov Spaces}

An interesting family of subspaces are the \definition{Krylov spaces}
\[
  \Krylov(\A, \v,n) = \Span(\v,\A \v,\A^2 \v,\ldots,\A^{n-1} \v)
\]
where \(\v\neq0\) is some arbitrary vector.  The Ritz values of \(\A\) in a
Krylov space are non-degenerate, as otherwise we could construct a Ritz vector
\(\y\neq0\) such that \((\y,\v) = 0\): this would imply that \((\y,\A^k\v) =
\theta^k(\y,\v)=0\;\forall k\geq0\) contradicting the assumption that
\(\y\in\Krylov(\A,\v,n)\).

A similar argument shows that for Ritz pairs \(\theta,\y\) and \(\theta',\y'\)
there is a unit vector \(\u\in\Span(\y,\y')\) with \((\u,\v) = 0\), for which
\((\u,\A^k\v)\propto(\theta-\theta')\|\v\|\), indicating that although we may
find good approximations in a Krylov space for nearly-degenerate eigenvalues,
we should not expect the same to be true for the corresponding eigenvectors.

\section{Error Bounds}

\subsection{Kaniel--Paige--Saad Bounds} \label{sec:KPS}

In this section we briefly survey the derivation of the Kaniel--Paige--Saad
bounds on the extremal eigenvalues of \Hermitean\ matrices.  Our notation is
similar but not identical to that of Parlett~\cite{Parlett:280490}.

Let us now study how rapidly we may expect eigenvalues to converge to a given
accuracy in the sequence of Krylov spaces of increasing dimension.  We shall
follow the method introduced by Kaniel~\cite{Kaniel:1966} and corrected and
extended by Paige~\cite{Paige:1971} and Saad~\cite{Saad:1980}.  Consider the
angle \(\angle(\z_j,\u)\) between an eigenvector \(\z_j\) and a vector
\(\u\in\Krylov(\A,\v,n)\).  We shall assume that the initial vector is not
orthogonal to the eigenvector, \((\z_j,\v)\neq0\): this would be very unlikely
for a randomly chosen \(\v\).  Indeed, we expect than on average \(|(\z_j,\v)|
\approx1/\sqrt N\) if \(\v\) is a random unit vector.

Consider the subspace \(X\subseteq\Krylov(\A,\v,n)\subseteq\R^N\), with
\(X=\Span(\y_j,\y_{j+1},\ldots,\y_n)\) where the \(\y_i\) are Ritz vectors for
\(\Krylov(\A-\lambda_j,\v,n) = \Krylov(\A,\v,n)\) for which the MaxMin bound
gives the first and last inequalities in
\[
  \theta_j-\lambda_j
    = \max_{S_j \subseteq X}\min_{\u \in S_j}\rho(\u,\A-\lambda_j)
    \leq \max_{S_j \subseteq\R^N}\min_{\u \in S_j}\rho(\u,\A-\lambda_j)=0,
\]
the central inequality following from the inclusion \(X\subseteq\R^N\).
\(\rho(\u\in X,\A)\leq\theta_j\) as it is a weighted mean of \(\theta_j,
\theta_{j+1},\ldots,\theta_n\), hence \(\rho(\u\in X,\A-\lambda_j)\leq
\theta_j-\lambda_j\), or equivalently
\[
  0\leq\lambda_j-\theta_j\leq\rho(\u\in X,\lambda_j-\A).
\]

Any \(\u\) in the Krylov space can be written as \(\u=p(\A)\v\) where \(p\) is
a polynomial of \(\deg p\leq n-1\).  A sufficient condition for \(\u=p(\A)
\v\) to be in \(X\) is that the polynomial \(p\) has the Ritz values
\(\theta_1,\theta_2,\ldots,\theta_{j-1}\) as roots, \(p(\theta_i)=0\) for
\(i<j\).  If \(p(\A)=(\A-\theta_i)q(\A)\) with \(q\) some polynomial of degree
\(\leq n-2\) then \(\bigl(\y_i,\u\bigr)=\bigl(\y_i,p(\A) \v\bigr)=\bigl(\y_i,
(\A-\theta_i)q(\A)\v\bigr)=\bigl((\A-\theta_i)\y_i, q(\A)\v\bigr)=(\w,q(\A)\v)
=0\), where \(\w\) is a vector orthogonal to \(\Krylov(\A,\v,n)\) and hence to
\(q(\A)\v\in \Krylov(\A,\v,n-1)\).  Therefore a sufficient condition for
\(\u\in X\) is that \(p(x)=q(x)l(x)\) with \(\deg q\leq n-j\) and \(l(x)=
\prod_{i=1}^{j-1}\frac{\theta_i-x}{\theta_i-\lambda_j}\).  Although we do not
need it here it is easy to prove that this is also a necessary condition.

The Rayleigh quotient for a vector \(\u=p(\A)\v=\sum_{k=1}^Np(\lambda_k)\z_k
(\z_k,\v)\) is
\[
  \rho(\u\in X,\lambda_j-\A)
    = \frac{\bigl(\u,(\lambda_j-\A)\u\bigr)}{(\u,\u)}
    = \frac{\sum_{k=1}^N |p(\lambda_k)(\z_k,\v)|^2(\lambda_j-\lambda_k)}
      {\sum_{k=1}^N |p(\lambda_k)(\z_k,\v)|^2}.
\]
The terms in the numerator for \(k<j\) are all negative, so we may bound
\(\rho\) above by dropping them; likewise we only increase the bound by
replacing \(\lambda_j-\lambda_k\) by \(\lambda_j-\lambda_N\) for the remaining
terms in the numerator, and dropping all the terms with \(k\neq j\) in the
denominator, leading to
\[
  \rho(\u\in X,\lambda_j-\A)
    \leq(\lambda_j-\lambda_N) \frac{\sum_{k={j+1}}^N
      |p(\lambda_k)(z_k,\v)|^2}{|p(\lambda_j)(\z_j, \v)|^2}.
\]
From this we immediately obtain the desired bound using Kaniel's choice of
\(q\) as the Chebyshev polynomial \(T_{n-j}\circ\gamma_j\), where
\[
  \gamma_j(\lambda)
    = 2\frac{\lambda-\lambda_N}{\lambda_{j+1}-\lambda_N}-1 
    = 1+2\frac{\lambda-\lambda_{j+1}}{\lambda_{j+1}-\lambda_N}
\] 
maps the interval \([\lambda_N,\lambda_{j+1}]\mapsto[-1,1]\), hence (q.v.,
Appendix~\ref{sec:app})
\[
  \max_{\lambda\in[\lambda_N,\lambda_{j+1}]}|q(\lambda)|
    =\max_{x\in[-1,1]}|T_{n-j}(x)| \leq1
\]
Furthermore \(l(\lambda_j)=1\), so if we define \(K_j\defn|l(\lambda_N)|\geq
\max_{\lambda\in[\lambda_N,\lambda_{j+1}]}|l(\lambda)|\) we obtain
\[
  \rho(\u\in X,\lambda_j-\A) 
    \leq\frac{(\lambda_j-\lambda_N)(K_j\tan\angle_j)^{2} }
             {T_{n-j}\big(\gamma_j(\lambda_j)\big)^2}.
\]
Here \((\tan\angle_j)^2=\frac{\sum_{k=j+1}^N|(\z_k,\v)|^2}{|(\z_j,\v)|^2}\),
so \(\angle_j\) is the angle between the projection of \(\v\) onto
\(\Span(\z_{j+1},\ldots,\z_N)\) and \(\z_j\).  Note that in the interesting
case where \(\theta_i\approx\lambda_i\) for \(i\leq j\)
\[
  K_j = \prod_{i=1}^{j-1}\frac{\theta_i-\lambda_N}{\theta_i-\lambda_j}
  \approx \prod_{i=1}^{j-1}\frac{\lambda_i-\lambda_N}{\lambda_i-\lambda_j}.
\]

From the bounds on Chebyshev polynomials given in Appendix~\ref{sec:app} we
get
\[
  0 \leq \lambda_j-\theta_j 
  \leq 4(\lambda_j-\lambda_N)(K_j\tan\angle_j)^2e^{-4(n-j)\mu_j+O(\mu_j^2)}
\]
where \(\gamma_j(\lambda_j)=1+2\mu_j^2\) with \(\mu_j=\sqrt{\frac{\lambda_j-
\lambda_{j+1}} {\lambda_{j+1}-\lambda_N}}\), and thus \(\ln\Bigl(\gamma_j
(\lambda_j) +\sqrt{\gamma_j(\lambda_j)^2-1}\Bigr)=2\mu_j+O(\mu_j^2)\).

\subsection{Bounds for Interior Eigenvalues} \label{sec:s-s-bounds}

We will produce a bound for an interior eigenvalue by introducing a shift
\(\Lambda\) and considering the spectrum of the shifted and squared matrix
\(\A'\equiv-(\A-\Lambda)^2\).  The parameter \(\Lambda\) is only needed for
our error bound and has no effect on the Krylov space Ritz values; hence the
best bound will just be the minimum over all possible values of \(\Lambda\); a
practical bound can be found by a judicious choice of \(\Lambda\).

Notice that the ``interior'' eigenvalues of \(\A\) near \(\Lambda\) are the
largest eigenvalues of \(\A'\).  Furthermore, the Krylov space for \(\A'\) is
contained within that for \(\A\), albeit with twice the size (but the same
number of applications of \(\A\)): \(\Krylov(-(\A-\Lambda)^2,\v,n) \subseteq
\Krylov(\A, \v,2n)\).  We will establish that the Kaniel--Paige--Saad bounds
for \(\Krylov(\A',\v,n)\) must therefore also hold for the latter.

Suppose that the spectrum of \(\A\) has an interior eigenvalue
\(\lambda_\alpha\).  As usual, we shall write the eigenvalues in decreasing
order \(\lambda_1\geq\cdots\geq\lambda_{\alpha-1}\geq\lambda_\alpha\geq\cdots
\geq\lambda_N\), as we shall do for the corresponding eigenvalues of \(\A'\),
\(\lambda'_1\geq\cdots\geq\lambda'_{j-1}\geq\lambda'_j=-(\lambda_\alpha-
\Lambda)^2\geq\lambda'_{j+1}\geq\lambda'_N\).  The smallest eigenvalue of
\(\A'\) is \(\lambda'_N=\min\bigl(-(\lambda_1-\Lambda)^2,-(\lambda_N-
\Lambda)^2\bigr)\).

We apply the Kaniel--Paige--Saad bound to the matrix \(\A'\) in the Krylov
space \(\Krylov(\A',\v,n)\) whose Ritz pairs are \(\theta'_j,\y'_j\) for
\(j=1,\ldots,n\).  Following the arguments of the preceding section we find
that the Ritz value \(\theta'_j\) satisfies \(\lambda'_j \geq\theta'_j\geq
\rho(\u\in X',\A')\) for any vector \(\u\in X'\defn\Span(\y'_j,\y'_{j+1},
\ldots,\y'_n)\) and the Rayleigh quotient is bounded by choosing \(\u=p'(\A')
\v\) for the polynomial
\[
  p'(x)=T_{n-j}\biggl(2\frac{x-\lambda'_N}{\lambda'_{j+1}-\lambda'_N}-1\biggr)
    \prod_{i=1}^{j-1}\frac{\theta'_i-x}{\theta'_i-\lambda'_j}.
\]
This leads to \(0\leq\lambda'_j-\theta'_j\leq\rho(\u\in X',\lambda'_j-\A')\)
with
\begin{equation}
  \rho(\u\in X',\lambda'_j-\A')
    \leq 4(\lambda'_j-\lambda'_N)(K'_j\tan\angle'_j)^2e^{-4(n-j)\mu'_j
      +O({\mu'_j}^2)}
    \defn\varepsilon'_{\Lambda},
  \label{eq:ss-bound}
\end{equation}
where
\[
  K'_j \defn \prod_{i=1}^{j-1}\frac{\theta'_i-\lambda'_N}{\theta'_i-\lambda'_j}
    \approx \prod_{i=1}^{j-1}\frac{\lambda'_i-\lambda'_N}
	    {\lambda'_i-\lambda'_j},\qquad
  \mu'_j=\sqrt{\frac{\lambda'_j-\lambda'_{j+1}} {\lambda'_{j+1}-\lambda'_N}},
\]
and \(\angle'_j\) is analogous to \(\angle_j\) (the eigenvectors of \(\A'\)
are the same as those of \(\A\) but their ordering is different).

Let \(\nu\defn|\lambda_\alpha-\Lambda|\), \(\xi\defn\sqrt{-\theta'_j}>0\)
(recall that with our conventions both \(\lambda'_j\leq0\) and \(\theta'_j\leq
0\)), and \(\delta\defn\lambda'_j-\theta'_j=-\nu^2+ \xi^2=(\xi-\nu)(\xi+\nu)\)
with \(0\leq\delta\leq\varepsilon'_\Lambda\).  Then \(\xi-\nu=\frac\delta{\xi+
\nu}\leq\frac\delta{2\nu}\leq\frac{\varepsilon'_\Lambda}{2\nu}\) since
\(\delta\geq0\implies\xi-\nu\geq0\implies\xi+\nu=(\xi-\nu)+2\nu\geq2 \nu\).
Hence we obtain the following bound on the error in the Ritz value
\(\bar\theta_\alpha=\Lambda\pm\sqrt{-\theta'_j}\) of \(\A\) corresponding to
the eigenvalue \(\lambda_\alpha\) of \(\A\) obtained from the Ritz value
\(\theta'_j\) of \(\A'\) in the Krylov space \(\Krylov(\A',\v,n)\):
\begin{equation}
  |\lambda_\alpha-\bar\theta_\alpha|
  = \min_\pm\left|\lambda_\alpha-\left(\Lambda\pm\sqrt{-\theta'_j}
    \right)\right|
  = \xi-\nu\leq\frac{\varepsilon'_\Lambda}{2\nu}
  = \frac{\varepsilon'_\Lambda}{2|\lambda_\alpha-\Lambda|}.
  \label{eq:bar-bound}
\end{equation}
Furthermore, since \(\Krylov(\A',\v,n)\subseteq\Krylov(\A,\v,2n)\subseteq
\R^N\) the MaxMin inequalities
\[
    \max_{S_\alpha\subseteq \Krylov(\A',\v,n)}\min_{\u\in S_\alpha}\rho(\u,\A)
    \leq\max_{S_\alpha\subseteq \Krylov(\A,\v,2n)}\min_{\u\in S_\alpha}
      \rho(\u,\A)
    \leq \max_{S_\alpha\subseteq\R^N}\min_{\u\in S_\alpha}\rho(\u,\A)
\]
imply that \(\bar\theta_\alpha\leq\theta_\alpha\leq\lambda_\alpha\).  This
means that the bound (\ref{eq:bar-bound}) on the error of the Ritz value
\(\bar\theta_\alpha\) from the Krylov space \(\Krylov(\A',\v,n)\) bounds that
of the Ritz value \(\theta_\alpha\) from the Krylov space \(\Krylov(\A,\v,
2n)\), \(\lambda_\alpha-\bar\theta_\alpha\geq\lambda_\alpha-\theta_\alpha
\geq0\).

\subsection{Discussion}

The advantage of the bound given by equations (\ref{eq:ss-bound})
and~(\ref{eq:bar-bound}) is that although the exponents \(\mu'\) will in
general be smaller than those from the original Kaniel--Paige--Saad bound the
prefactor \(K'\) may be far smaller in the interior of the spectrum.

Furthermore, \(\mu'\) may not be much smaller than \(\mu\) for eigenvalues in
low density parts of the spectrum for a good choice of the parameter
\(\Lambda\).  For example, suppose that \(\lambda_\alpha\) is is just below a
large void in the spectrum of width \(\Delta\), but not necessarily
immediately adjacent to it, and is also close to \(\lambda_{\alpha+1}\).  We
can choose \(\Lambda\) near the middle of void such that \(\lambda'_{j+1}=
-(\lambda_{\alpha+1}-\Lambda)^2\).  For notational simplicity let us also
assume that \(\lambda'_N=-(\lambda_N-\Lambda)^2\), so
\[
  {\mu'_j}^2 = \frac{\lambda'_j-\lambda'_{j+1}} {\lambda'_{j+1}-\lambda'_N}
  = \frac{-(\lambda_\alpha-\Lambda)^2+(\lambda_{\alpha+1}-\Lambda)^2}
      {-(\lambda_{\alpha+1}-\Lambda)^2+(\lambda_N-\Lambda)^2}
  = \frac{\bigl(\lambda_\alpha-\lambda_{\alpha+1}\bigr)
      \bigl(2\Lambda-\lambda_\alpha-\lambda_{\alpha+1}\bigr)}
    {\bigl(\lambda_{\alpha+1}-\lambda_N\bigr)
      \bigl(2\Lambda-\lambda_{\alpha+1}-\lambda_N\bigr)}.
\]
Now, \(2\Lambda-(\lambda_\alpha+\lambda_{\alpha+1})\approx\Delta\), hence
\[
  \mu'_j = \mu_j \sqrt{\frac{2\Lambda-\lambda_\alpha-\lambda_{\alpha+1}}
    {2\Lambda-\lambda_{\alpha+1}-\lambda_N}}
  \approx \mu_j \sqrt{\frac{\Delta}{2\Lambda-(\lambda_{\alpha+1}+\lambda_N)}}
  \geq \mu_j \sqrt{\frac\Delta{\lambda_1-\lambda_N}},
\]
which is only slightly smaller than \(\mu_j\) if \(\Delta\) is a significant
fraction of the spectrum.

The limit on how far \(\Lambda\) can be from \(\lambda_\alpha\) and
\(\lambda_{\alpha+1}\) is that when it gets closer to a dense band of
eigenvalues the opposite side of the void the ordering of the eigenvalues of
\(\A'\) will change and \(\lambda'_{j+1}<-(\lambda_{\alpha+1}-\Lambda)^2\).  A
few eigenvalues within the void do not matter, as they can be ``removed'' by a
suitable zero in the polynomial \(p'\), at the cost of a smaller \(K'\) factor
and a slightly smaller factor \(n-j\) in the exponent.

Amusingly, for (near) extremal eigenvalues \(\mu=\mu'\) as we may
take~\(\Lambda\to\pm\infty\).

Nevertheless, there is always an extra factor of two in the exponent arising
from the fact the Krylov space \(\Krylov(\A,\v,2n)\) is twice as large as
\(\Krylov(\A',\v,n)\), so asymptotically for large \(n\) the original bound is
always better.  This is not too important for two reasons: firstly because we
are usually interested in getting approximations within some specified error
rather than in the limit where the error tends to zero, and secondly since the
Ritz values are exactly equal to the eigenvalues for \(n=N\) in exact
arithmetic (or for smaller values of \(n\) if the spectrum of \(\A\) is
degenerate) there is a natural upper limit on the dimension \(n\) of the
Krylov space.

\section{Examples}

\subsection{Qualitative Example}

As our first example, consider an \(N\times N\) matrix \(\A\) whose spectrum
lies in the unit interval.  This matrix has \(k\) bands, each containing
\(N/k\) eigenvalues, and each band is separated by a ``void'' of
size~\(\Delta\).  Within each band the eigenvalues are separated uniformly by
\(\delta=[1-(k-1)\Delta]/N\).

If \(\Delta=0\) then all the eigenvalues are uniformly spaced, so \(\mu_j
\propto\sqrt\delta=1/\sqrt N\) for any \(\lambda_j\), and the
Kaniel--Paige--Saad bound shows that to resolve the eigenvalues we require a
Krylov space of dimension \(n\propto\sqrt N\).  For eigenvalues in the
interior of the spectrum the \(K_j\) factors become large, of course.

If \(\delta=0\) then \(\A\) has \(k\) eigenvalues each with a degeneracy of
\(N/k\), and there are only \(k\) Ritz pairs.  If \(\delta\ll1\) then the
error in the Ritz values \(|\theta-\lambda|\approx\delta\) for \(n\approx k\),
but we require \(n\propto1/\sqrt\delta\) in order to resolve the eigenvalues
within each band (at least until the Ritz values become exact for~\(n=N\)).

If \(\Delta\gg\delta\) then the exponent \(\mu'\approx\sqrt\delta\) for an
interior eigenvalue \(\lambda\) close to a void for shift \(\Lambda\approx
\lambda\pm\Delta/k\), so to resolve it we also need \(n\propto\sqrt N\).  This
is the same scaling behaviour as for the usual bound, but the key difference
between our bound and the usual one is that~\(K'\ll K\).  As
\(\Delta\to\delta\) the scaling behaviour crosses over to \(n\propto N\) since
\(\mu'\approx\mu^2\).  Perhaps this should not be too surprising as the Ritz
values become exact for~\(n=N\).

\subsection{Quantitative Example}

\begin{figure}[ht!]
  \begin{center}
    \epspdffile{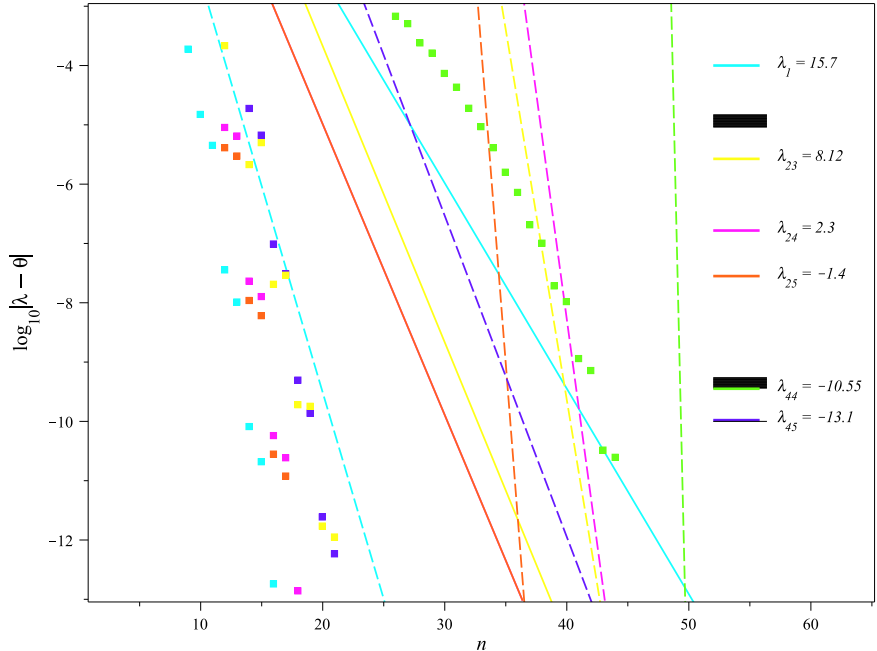}
  \end{center}
  \caption{Convergence of Ritz values for the matrix \(\A\) whose spectrum is
    shown on the right.  The plot shows the logarithm of the magnitude of the
    difference between each Ritz value and the nearest eigenvalue of \(\A\) as
    a function of the dimension \(n\) of the Krylov space.  The initial vector
    for the Krylov spaces was taken to have exactly equal overlap with each
    eigenvector.  The black ``bands'' in the spectrum are actually
    equally-spaced eigenvalues with a spacing of~0.05.  The smallest
    eigenvalue is \(\lambda_{46}=-13.2\).  The squares show the errors for the
    actual Ritz values, and their colour indicates the nearest eigenvalue.
    The solid lines are the bounds from the shifted and squared matrix \(\A'\)
    with shifts \(\Lambda\) chosen to give the smallest bound for an error of
    \(10^{-8}\), while the dashed lines are the original Kaniel--Paige--Saad
    bounds.  The shifts used were \(\Lambda\to\infty\) for \(\lambda_1\) and
    \(\Lambda=0.45\) for \(\lambda_{23}\), \(\lambda_{24}\),
    and~\(\lambda_{25}\).  The solid line for \(\lambda_{24}\) lies under that
    for \(\lambda_{25}\): this is because they are equidistant from the shift
    value.  The Ritz values coincide with the eigenvalues for
    \(n=\dim\A=46\).}
  \label{fig:bounds}
\end{figure}

Our second example is more quantitative.  We consider a matrix with two bands
of closely-spaced eigenvalues and a few other eigenvalues well-separated from
these bands.  In order to determine the value of the shift \(\Lambda\) we need
to specify the value of \(n\) --- the dimension of the Krylov space
\(\Krylov(\A,\v,n)\) --- at which the bounds are to be evaluated.  We do this
implicitly by choosing the value of \(\Lambda\) for which a bound of some
specified magnitude, in this case \(10^{-8}\), is achieved for the
smallest~\(n\).  We compute an approximation to such a \(\Lambda\) by
minimizing the exponent \(\mu'\).  To do this we note that since
\(d\mu'/d\Lambda\neq0\) anywhere \(\mu'\) can only have an extremum as a
function of \(\Lambda\) when two eigenvalue of \(\A'\) become degenerate, as
at such points the correspondence between the eigenvalues of \(\A\) and those
of \(\A'\) changes so as to maintain the correct order of the~\(\lambda'\).
We stress that these approximations merely mean that we may not have chosen
quite the optimal value for \(\Lambda\), nevertheless the bounds obtained are
still valid since they hold for any~\(\Lambda\).

The convergence of the Ritz values is illustrated in Figure~\ref{fig:bounds}.
This clearly shows that our new bounds (solid lines) are significantly tighter
than the usual Kaniel--Paige--Saad bounds (dashed lines) for the
well-separated interior eigenvalues, although as expected the original bounds
are better for the extremal eigenvalues.

This example illustrates some other interesting properties of the bounds.
Observe that the lowest eigenvalue in this example, \(\lambda_{46}=-13.2\), is
almost degenerate with the second-lowest eigenvalue, \(\lambda_{45}=-13.1\).
This means that its bounds correctly indicate that it will only converge
slowly.  Interestingly the bound for latter is better than that for the former
because this small gap appears in the coefficient \(K_2\) rather than the
exponent \(\mu_1\) of the extremal bounds for \(-\A\).

If we did not know the spectrum beforehand, but tried to determine the bounds
from the Ritz values from \(\Krylov(\A,\v,n)\) for some large value of \(n\),
we might be misled into thinking that the lowest eigenvalue was exactly
degenerate, and thus obtain a much tighter bound on its convergence, and this
bound would seem to agree with the empirical convergence of the Ritz value.
Only when \(n\) becomes sufficiently large to resolve the splitting is this
bound violated, there being a ``plateau'' in the error for this Ritz value
until the Krylov space is able to separate the actual eigenpairs.  The fact
that the eigenvalue is nearly degenerate may be determined even for small
\(n\) by ``restarting'', that is constructing the Krylov space
\(\Krylov(\A,\v',n)\) for some different initial vector \(\v'\), but this is
insufficient to distinguish between exact and near degeneracies.

It is often the case in applications that we are interested in finding
eigenvalues and their eigenspaces up to some specified accuracy: for example,
the matrix may be derived from a model that is only valid up to this accuracy.
In such cases the bounds determined from the approximate spectrum may be more
useful than the ``exact'' bounds, as might the ``exact'' bounds obtained from
the approximate spectrum in which nearly degenerate eigenvalues are collapsed
into exactly degenerate ones.  This situation is common to both our shifted
and squared bounds and to the original ``extremal'' bounds.

\section{Acknowledgements}

We gratefully acknowledge the support of the Centre for Numerical Algorithms
and Intelligent Software (EPSRC EP/G036136/1) in the preparation of this work.

\appendix\section{Chebyshev Polynomials} \label{sec:app}

\definition{Chebyshev polynomials} of the first kind may be defined as
\(T_n(x)=\cos(n\cos^{-1} x)\), so it is immediately obvious that \(|T_n(x)|
\leq1\) for \(|x|\leq1\).  It is perhaps surprising that it is a polynomial in
\(x\), but application of the binomial theorem shows that
\[
  T_n(x) = \sum_{j=0}^{\lfloor n/2\rfloor}\sum_{\ell=0}^j \binomial(2n,j)
    \binomial(j,\ell)(-1)^{j+\ell}x^{n+2(\ell-j)},
\] 
so it is indeed a polynomial of degree~\(n\).  Moreover it is also easy to see
that
\[
  T_n(x) = \half\left[ e^{n\ln\left(x+\sqrt{x^2-1}\right)}
    + e^{n\ln\left(x-\sqrt{x^2-1}\right)}\right];
\]
we thus deduce that \(T_n(x)\geq\half e^{n\ln\left(x+\sqrt{x^2-1}\right)}\)
for \(x\geq1\).
\bibliographystyle{siam}
\bibliography{lanczos-theory}
\end{document}